\documentclass[12pt]{amsart}

\newcommand{\invisible}[1]{}

\newcommand{\numbering}{\noindent {\addtocounter{subsection}{1}{\bf \thesubsection}}---}

\author[C. Simpson]{Carlos Simpson}
\address{CNRS, Laboratoire J. A. Dieudonn\'e\\ Universit\'e de Nice-Sophia Antipolis\\
06108 Nice, Cedex 2, France}
\email{carlos@math.unice.fr}
\urladdr{http://math.unice.fr/$\sim$carlos/} 
\title[Theorem proving]{Computer theorem proving in math}

\begin{document}

\begin{abstract}
We give an overview of issues surrounding computer-verified 
theorem proving in the standard pure-mathematical context. 
\end{abstract}

\keywords{Lambda calculus, Type theory, Theorem proving, Verification, Set theory}

\maketitle

\section{Introduction}
\label{introduction}

When I was taking Wilfried Schmid's class on variations of
Hodge structure, he came in one day and said ``ok, today
we're going to calculate the sign of the curvature
of the classifying space in the horizontal directions''.
This is of course the key point in the whole theory: the
negative curvature in these directions leads to all sorts
of important things such as the distance decreasing property.

So Wilfried started out the calculation, and when he came
to the end of the first double-blackboard, he took the
answer at the bottom right and recopied it at the upper
left, then stood back and said ``lets verify what's
written down before we erase it''. Verification made (and
eventually, sign changed) he erased the board and started
in anew. Four or five double blackboards later, we got to
the answer. It was negative. 

Proof is the fundamental concept underlying mathematical research. 
In the exploratory mode, it is the main tool by which we percieve the mathematical world as it {\em  really is}
rather than as we would {\em like it to be}. 
Inventively, proof is used for validation of ideas: one
uses them to prove something nontrivial, which is valid only if the proof is
correct.  Other methods of validation exist, for example showing that an idea leads to
computational predictions---but they often generate proof obligations too. 
Unfortunately, the need to prove things is the factor which slows us down the most too. 

It has recently become possible, and also necessary, to imagine a full-fledged 
machine-verification 
system for mathematical proof. This might radically change many aspects of mathematical research---for better
or for worse is a matter of opinion. At such a juncture it is crucial that standard pure mathematicians 
participate.

The present paper is intended to supply a synthetic contribution to the subject.  Most, or probably all,
is not anything new (and I take this opportunity to apologize in advance for any additional references which should
be included).  What is supposed to be useful is the process of identifying a certain collection of problems, and
some suggestions for solutions, with a certain philosophy in mind. 
The diverse collection of workers in this field provides a diverse collection of philosophical perspectives,
driving our perception of problems as well as the choice of solutions. Thus
it seems like a valid and useful task to set out how things look from a given philosophical point of view. 
Our perspective 
can be summed up by saying
that we would like to formalize, as quickly and easily as possible, the largest amount of standard mathematics,
with the ultimate goal of getting to a point where it is concievable to formalize current research mathematics
in any of the standard fields.  

This is a good moment to point out that one should not wish that everybody share the same point 
of view, or even anything close.  On the contrary, it is good to have the widest possible variety of questions
and answers, and this is only obtained by starting with the widest possible variety of philosophies.

We can now discuss the difference between what we might call ``standard'' mathematical practice, and
other currents which might be diversely labelled ``intuitionist'', ``constructivist'' or ``non-standard''.
In the standard practice, mathematicians feel free to use whatever system of axioms they happen to have
learned about, and which they think is currently not known to be contradictory. This can even involve
mixing and matching from among several axiom systems, or employing reasoning whose axiomatic basis is not
fully clear but which the mathematician feels could be fit into one or another axiomatic system if 
necessary.  This practice must be viewed in light of the role of proof as validation of ideas: for
the standard mathematician the ideas in question have little, if anything, to do with logical foundations,
and the mathematician seeks proof results for validation---it is clear that any generally
accepted framework will be adequate for this task.  

A growing number of people are 
interested in doing mathematics within an intuitionist or constructive framework. They
tend to be closest to mathematical logic and computer programming, which may be why many current
computer-proving tools are to some degree explicitly designed with these concerns in mind.
There are many motivations for this, the principal one being deeply philosophical.
Another is just concern about consistency---you never know if somebody might come up with an inconsistency
in any given axiomatic system (see \cite{Kitoda} for example). 
Common sense suggests that 
if you take a system with less axioms, there is less chance of it,
so it is possible to feel ``safer'' doing mathematics with not too many axioms. 
A subtler motivation is the question of knowing what is the minimal axiomatic system under
which a given result can be proven, see \cite{ReverseSimpson} for example.
Computer-scientists have some wonderful toys which make it possible 
directly to transform a constructive proof into a computational algorithm.

The problem of computer-verified-proof in the standard framework
is a distinct major goal.  Of course, if you prove something in a constructive framework,
you have also proven it in the standard framework but the opposite is not the case.  
Integrating additional axioms such as replacement or choice
into the theory, gives rise to a distinct optimization problem: it is often easier to prove a given thing
using the additional axioms, which means that the structure of the theory  
may be different. The notion of computer verification has already made important
contributions in many areas outside of ``standard'' mathematics, and it seems reasonable to think that
it should also provide an important tool in the standard world.

One could also note that, 
when you get down to the nitty-gritty details, most technical results concerning 
lambda calculus, including semantics, 
normalization and consistency results for various (highly intuitionistic) axiomatic systems,
are themselves proven in a standard mathematical framework relying on Zermelo-Frankel set theory 
(see \cite{WernerThesis} for
example). Thus even
people who are interested in the intuitionist or constructive side of things might find it useful to 
have good standard tools at their disposal. 

One reason why computer theorem verification in the standard pure mathematical context is not
recieving enough attention is that most pure mathematicians are unaware of the subject. This lack of
awareness---which I shared a few years ago---is truly colossal, given the literally {\em thousands of papers}
concerning the subject which have appeared in recent years. Thus it seems a bit silly, but one of the
goals of the present paper is to try to spread the news.

This ongoing work is being done with a number of other
researchers at Nice: A. Hirschowitz, M. Maggesi, L. Chicli,
L. Pottier. I would like to thank them for many stimulating interactions. 
The referee kindly pointed out the references \cite{Bundy} and \cite{KaufmannMoore}, 
and several readers of the first version sent helpful
remarks and references. I would like to thank
C. Raffalli, creator of the PhoX system, who explained some of its ideas during his visit to Nice
last summer---these have been very helpful in understanding how to approach proof documents. 
Special thanks are extended to the several members of the Coq group at INRIA, Rocqencourt, for their
continuing help during visits and on the mailing list. From much
further past, I would like to thank M. Larsen for teaching me what little I know about computer programming. 

I would like warmly to thank S. Gutt for organizing a great conference, for letting me give a talk which
was fairly widely off-subject, and for suggesting a title. 
And of course the birthday-teenagers for reminding us of the Sumerian numbering system \cite{Melville}. 

\section{Basic reasons}
The basic goal of mechanized theorem verification is to have a language for creating 
mathematical documents, such that the
reasoning can be checked by computer.\footnote{We ask that the machine verify the proofs supplied by the
mathematician; obviously it can help if the machine finds the proof on its own but that constitutes a
difficult, and somewhat distinct, research problem \cite{Bundy}.}
This is of interest to pure mathematicians for reasons such as the following. 

\numbering
To reduce the number of mistakes we make in our mathematical work (there are lots of them! \cite{Neeman}).

\numbering 
Dealing with complicated theories---it is difficult to manage
development when there are too many things to prove
(and often the author doesn't even know exactly what it is
he needs to check).

\numbering 
Because of (2.1) and (2.2), nobody really has the time or energy
adequately to verify all new results, which leads to
problems for publication or other timely dissemination. If articles could be validated
by machine then referees could concentrate on the question
of how interesting they are.

\numbering 
Black-box theories: it is sometimes useful, but
currently dangerous, to use results without knowing how they
are proven. An author can have a certain situation in mind, and might use arguments which work in
that situation but not in other ones, but sometimes the explicit statements of results don't adequately
reflect all the hypotheses. So future users need to understand the arguments before 
blindly applying the statements. An advantage of machine verification is that the language
has to be totally precise: utilization
of a ``black box''  via cut-and-paste rather than axiomatically 
leaves no room for misinterpretation.

\numbering 
Integration of many theories or points of view: it is hard
for humans because of the steep learning curve for each
new theory.
Even now the importation of {\em techniques},
without a full understanding of how they work, is a little-heralded part of mathematical practice,
generally glossed over because of the problems of accuracy that are posed. 
With machine-verified proving these problems should in large part be resolved and should lead 
to a much more diverse fabric in which
different workers develop different techniques and import the techniques of others without 
needing to spend such a long time learning about the details, leaving everyone with more energy to spend 
looking for new ideas.

\numbering 
New forms of reasoning: now, calculations or proofs with
little motivation are ruled out by the fact that the
mathematician needs to be motivated in order to understand
what is going on. Locally unmotivated reasoning 
might lead to globally interesting results. A striking example in logic is work on minimal
length axioms such as \cite{McCuneVeroff} \cite{XCB}, where much of the reasoning was found
using theorem-proving software.

\numbering 
The advent of more and more proofs making use of computer calculation (non-existence of a projective plane
of order $10$ \cite{LamThielSwiercz}
\cite{LamThielSwierczMcKay}, the
4-color theorem, more recently a theorem on hyperbolicity of $3$-manifolds
\cite{GabaiMeyerhoffThurston} and the Kepler conjecture \cite{Hales}) represents {\em a priori} a 
somewhat different use of the computer than what we are considering here. However, these proofs pose
acutely the problem of verifying logical details, and a formalized framework promises to present
a significant improvement. This issue was addressed by B. Werner in a talk at Nice about the 4-color
theorem, and has motivated Hales' ``flyspeck project'' \cite{Hales}.

\numbering 
The possibility of computer-assisted proof verification would significantly extend the 
age period in which mathematicians are optimally productive: older mathematicians tend to avert getting
bogged down in details, which often prevents full development of the continuing stream of valid ideas.

On the whole, it is a mistake to think that computer theorem-proving will always lag behind
regular math. At some point, the advantages will so much more than counterbalance the difficulties
that mathematicians using these tools will jump out ahead.

\section{History}

\numbering Cantor, Hilbert, Russell; then G\"{o}del and so on,
worked on providing a strict logical foundation for
mathematics.

 Church, Rosser, Martin-L\"{o}f \cite{MartinLof}, Girard \cite{Girard} worked more
specifically on the logical foundations of type theory.

\numbering Automath, the project of de Bruijn, was one of the first
explicit projects in this direction; notably he pointed
out how to deal with variables (via de Bruijn indices).

In 1977,
Jutting entered all of Landau's analysis book in the AUTOMATH system \cite{Jutting}.

\numbering The first long-term project involving a very large number of
people was Mizar; it is still ongoing, and led to
the {\em Journal of Formalized
Mathematics} which has been publishing articles verified in the MIZAR system
since at least 1989. This includes much basic theory but also research articles.
The mathematical orientation of the
research articles is sometimes distinguished from the mainstream of
what we call standard mathematics (this is an observation but not a criticism).

\numbering More recently, a number of systems based on some kind of
type theory have been developped: HOL, Lego, Isabelle,
Nuprl, Nqthm, AC2L, and general ``Boyer-Moore'' theorem
provers...

 Coq is a project which seems relatively successful with
respect to the considerations discussed below.

An overview of very recent systems would include: Elf,
Plastic, Phox, PVS, IMPS, Ghilbert, Metamath.  These and many more can be
found using an internet search engine---a little perseverance turns up an astounding amount of
reference material which we couldn't begin to include here (and which is impossible
to read in its entirety). 
A good place to start is with meta-pages such as \cite{OxfordPage}, \cite{pfenning}.

\numbering Industrial demand, spurred in part by some well-known 
problems that might have been avoided with use of the
appropriate proof technology (Pentium, Ariane 5), continues
to be a major  impetus for the development of proof engines.
In this case, interest is oriented toward proving the
correctness of algorithms, integrated circuits or other
wiring diagrams, and so forth. 

\numbering
More and more conferences, journals, and mailing lists
bear witness to the fact that fields related to computer proof
verification are really booming. Browsing the archives of mailing lists
such as \cite{qed} or \cite{coqclub} is a fascinating homework assignment. 

\numbering
In contrast we might mention that the notion of computer
calcuation systems, which is perhaps much better known in
the world of pure mathematics, is not quite the same thing
as computer proof verification.  It is surely related,
though, and is also experiencing fast growth. {\em Idem} for the
closely related field of {\em automated deduction} 
(see \cite{Bundy} for a survey).

\numbering Some recent examples of proof developments which are of interest to the
pure mathematician are: the fundamental theorem of algebra, by Geuvers {\em et al} \cite{fta};
R. O'Connor's formalization of G\"{o}del's theory \cite{OConnor}; 
the \verb}Reals} and other recent libraries for Coq \cite{Coq};  L. Chicli's thesis
about algebraic geometry \cite{Chicli}; the logical libraries
in Isabelle \cite{Isabelle}; Metamath's work on Hilbert spaces \cite{Metamath}, \ldots .

A glance at this history leads to the conclusion that it
now seems like all of the ingredients are available so that
we could start doing computer theorem verification in a serious
way. We just have to find the right combination.
It doesn't seem unreasonable to think that this moment will
be seen, in retrospect, as an important juncture in
the history of mathematics.

\section{Extra-mathematical considerations}

It is useful to elucidate some criteria explaining
what would make the difference between a system which would or would not
be used and accepted by mathematicians. The model 
of mathematical typesetting and \TeX suggests some criteria. 

\numbering {\sc free}: the program shouldn't cost anything (and should
be free of any other constraints) for
scientific users. 

\numbering {\sc plain text}: the main format of the proof document
should be a plain text file; this is to insure that
nothing is hidden to the user between the document creation
and its compilation or verification, and also to insure 
easy exchange of documents.

\numbering {\sc open source}: the basic programming code used to
specify the syntax of the language, and to verify the
proof document, should be available to
everybody. In the best case it would actually be attached
as a trailer to the proof document itself, since the
program is an integral part of the proof specification
(see the discussion of archiving in the introduction of
\cite{GabaiMeyerhoffThurston}). 
Open source for this key part of the system brings the
guarantee that your document will be readable or compilable
at any time in the future.  

\numbering {\sc reference manual}:  it should be possible to get started 
using only
the reference manual, and at the same time the manual should
answer even the most arcane questions (first-time users
tend to stumble upon arcane questions more often than one
might expect, cf \cite{coqclub}).  

\numbering {\sc cross-platform}: it should be {\em very easy}
to install and use on all platforms. One
of the most difficult aspects for a mathematician getting
started in this subject is just to get a system up and
running usably on whatever computer is available.
The average mathematician should be able to install the
system and type out a first text file which is then
successfully verified, proving a simple mathematical lemma 
such as uniqueness of the identity element in a monoid, 
in no more than a couple of hours. 

\numbering {\sc modifiable}:
the system should be modifiable by the user, in particular the user should choose (or invent, or re-invent)
his own axiomatic framework,
his notational system, his proof tactics in a fully programmable language,
and so on---these things shouldn't be hard-wired in.  In William Gibson lingo, we want something
that would give the ``Turing police'' nightmares. This information would be like the 
notion of a macro header to a TeX file.  People could exchange and copy their header material, but could
also modify it or redo it entirely.  This is very important for enabling creative
research.

\numbering {\sc encourage mathematicians}: the people who maintain
the system should encourage standard mathematicians to
write math in the system. This includes 
maintaining an archive, and 
publicizing the existence of the system within the standard pure
mathematical community.   

\numbering {\sc attentiveness to user suggestions}: 
pure mathematicians are the scientists
who work the most with proof, and it seems logical that 
their viewpoints would be important to development of a
fully functional system.  Input from pure mathematicians
is likely to be useful even in more far-flung areas of the
subject such as industrial applications. This could be seen
as compensation for the requirements that the system be free
and open-source. 

\section{Starting tasks}

The hardest part of the subject is the first basic
task which needs to be addressed: you need to form an
opinion about what you are going to want your documents
to look like \cite{Wenzel}.

This task is of course first and foremost treated by the
system designers. However, even within the framework of
a given system, the mathematician has to address this task
too, in determining the style he will use.

It is subject to a number of constraints and hazards.

\numbering 
Everything you write down has to have an absolutely well-defined
and machine-determinable meaning.

\numbering 
One must write things in a reasonably economical way,
in particular we need a reasonable system of notation for
the standard mathematical objects which are going to be
manipulated. 

\numbering 
It is better if any difficulties which are encountered,
have actual mathematical significance. An example of
a difficulty which doesn't seem to have much real mathematical
significance arises in type theory when we have two
types \verb+A+ and \verb+B+, with an element \verb+a:A+ and a
known (or supposed) equality \verb+A=B+. It often (depending on
the type system being used) requires additional
work in order to be able to say that \verb+a:B+. To the 
extent possible, the pure mathematician will want to choose
a format in which any such hurdle is seamlessly cured.  

\numbering 
At least a part of what is written should be comprehensible,
preferably upon simple inspection of the pure text file.

\numbering 
{\sc potential problem}: it is quite possible 
(and maybe even generally the case) to invent notation
which suggests a meaning different from the actual meaning.
When you prove a given statement, did you really prove what
you thought you did? 

Unfortunately there is no way for
the computer to verify this part. The conclusion is that
we require human comprehension of the definitions, and
of the statements of the theorems (the main ones, at 
least). This imperious requirement is stronger than
the natural desire also to understand the proofs. 

\numbering 
Experience suggests that the proofs, in their explicit and
verifiable form, actually contain too much information and
the reader will not really want completely to understand 
them. Indeed the origin of the whole project is the fact
that we aren't able to follow the details of our
proofs. Rather one might ask that two particular
aspects stand out: 
the general strategy, and the
unexpected main steps.

\section{How it works}

We describe in brief some of the salient points of
the \verb+Coq v8.0+ system \cite{Coq} currently used by the author.
A number of other recent systems are similar.  Aside from
the geographical reason for my choice of system,
one can point out that \verb+Coq+ scores pretty well (but naturally
with room for improvement) on the above lists of 
{\em desiderata}. 
\footnote{If the makers of other systems feel 
impelled at this point to say ``hey, our system does just
as well or better'', the main point I am trying to 
make is: That's great! Try to get us working on it!}

\numbering
Some natural considerations lead to the idea of using
{\em type theory} as a basis for mathematical verification \cite{MartinLof}.
To start with, notice that we want to manipulate various
sorts of objects, first and foremost the mathematical 
statements themselves. Although perhaps not strictly necessary,
it is also very important to be able directly to manipulate 
the mathematical entities appearing in the statements, be
they numbers, sets, or other things. 
To deal with our entities (which will be denoted by
expressions entered into the text file), 
we need function application \verb+(f x)+ 
and function-definition
\begin{verbatim}
fun x => (... x ... x ....).
\end{verbatim}
This combination
of operators yields the {\em lambda calculus}. The
notion of {\em $\beta$-reduction} is that if
\begin{verbatim}
f:= fun x =>(... x ... x ....)
\end{verbatim}
then for any expression \verb+y+, the application \verb+(f y)+ should
reduce to 
\newline
\verb+(... y ... y ....)+. It is natural to ask that the computer take into
account this reduction automatically. 

\numbering
The first famous difficulty is
that if \verb}f:= fun x =>x x}
then \linebreak
\verb+(f f)= (fun x =>x x) f+ $\stackrel{\beta}{\rightarrow}$
\verb+(f f)+ $\stackrel{\beta}{\rightarrow}$ \verb+(f f)+ \ldots .
Thus the $\beta$-reduction operation can loop.

The solution  is the notion
of {\em type} introduced by Russell: every ``term'' (i.e. allowable
expression, even a hypothetical or variable one) has a
``type''. The meta-relation that a term \verb}x} has type \verb}X} is
written \verb}x:X}. The $\lambda$ construction (construction of
functions) is abstraction over a given type: we have
to write 
\begin{verbatim}
f:= fun (x:A)=>(... x ... x ....)
\end{verbatim}
where \verb+A+ is supposed to designate a type. Furthermore
the expression on the right is required to be allowable,
i.e. well-typed, under the constraint that \verb+x+ has type \verb+A+.

\numbering
{\sc The Church-Rosser theorem} states that:
\newline
---{\em In {\em typed} $\lambda$-calculus, $\beta$-reduction 
terminates; also calculation of the typing relation
terminates.}

The verification program (such as Coq) basically implements this calculation. 
The fact that it is sure to terminate 
is a considerable advantage for the end-user: you are sure
that the computer will stop and say either {\em ok} or {\em no}. And if ``no''
you were the one who made the mistake. 

\numbering
{\bf The product operation $\Pi$}: in the above expression,
what is the type of \verb+f+? We need a new operator $\Pi$ written as
\verb}f:forall (x:A), B x}
where 
{\small
\begin{verbatim}
B:=fun (x:A) => "the type of (...x ... x ....) when x has type A"
\end{verbatim}
}
(recall that the allowability restriction on \verb+f+
requires that this expression \verb+(B x)+ exist, and in fact
it is calculable). 

\numbering
{\bf Sorts:} What is the type of a type? 
We are manipulating the types themselves
as terms, so we are required to specify the type of a type.
This is known as a ``sort'' (or in some literature, a 
``kind'').
In current theories there are just a few, for example:
\newline
\verb+Prop+, the type of propositions or statements; and
\newline
\verb}Type_i}, the type of types which are thought of as sets
(with a ``universe index'' \verb+i+, see below), so an \verb}X:Type_i}
corresponds to a collection of objects, and
an \verb+x:X+ corresponds to an element of the collection.

\numbering
{\bf Curry-Howard:} In many (but not all) type-theory based systems,
the notion of mathematical proof comes out for free from the type-checking
algorithm. This is known as the ``Curry-Howard principle'': propositions are
realized as types, which the reader can think of as being sets with either zero
or one element. Proving a proposition amounts to giving an element of the corresponding
type, so the proof-checking mechanism is the same as the mechanism implementing the
Church-Rosser algorithm refered to above, determining whether a given expression is indeed
a term in the propositional type claimed by the mathematician. As suggested by the notation,
the product operation $\Pi$ corresponds to the logical ``forall'' in propositional calculus
(a product of zero-or-one-element sets is nonempty if and only if
all of the component sets are nonempty). 

The above system corresponds to Martin-L\"{o}f's type theory
\cite{MartinLof}. 
One can imagine theories with a more complicated diagram
of sorts, which is probably an interesting research topic
for computer scientists.

\numbering 
{\sc inductive objects:}
An interesting feature of a number of systems including
 \verb+Coq+ is the notion of inductive objects, which
allows us automatically to create recursive datatypes.
Since understanding and calculating with recursive datatypes
is what the computer does best, it can be good to make
use of this feature.  It is often the only viable way 
of getting the computer hitched onto the calculational
aspects of a problem. (This might change, notably with
the advent of connectivity between theorem-proving systems and 
calculation systems.)
We won't go into any further detail on the basic structure
or functionality of inductive definitions; the reader is
referred to the reference manual of \cite{Coq}.  

\section{Universes}

In order to get a coherent system (i.e. one without any
known proofs of \verb+False : Prop+) the sorts \verb}Type_i} have
to be distinguished by universe levels \verb+i+, with
\verb}Type_i : Type_(i+1)}
whereas \verb}Type_i} $\subset$ \verb}Type_(i+1)} also.
In the set-theoretical semantics,
these universe levels can be thought of as being Grothendieck universes,
see \cite{Werner}. 
The incoherency of the declaration \verb+Type:Type+ without
universe indices was shown in \cite{Girard}.  Nonetheless,
in some type-theoretical programs such as \verb+Coq+, a 
vestige of \verb+Type:Type+ is preserved by a neat, but quirky,
mechanism known as {\em typical ambiguity} (see \cite{Feferman}).
The user just writes \verb+Type+
and the computer makes sure that it is possible to assign indices
to all occurences of this word in such a way that
the required inequalities hold.  If this is not possible,
then a \verb+Universe Inconsistency+ error is raised. Thus, the
user {\em writes his text file as if} he were using
the inconsistent \verb+Type:Type+ system, and the verification
program then verifies that no unallowable advantage was
taken.  

Curiously enough, this leads to a phenomenon that looks
for all the world to me like quantum 
mechanics.\footnote{This might be related to \cite{Schmidhuber}.} 
Suppose you
prove \verb+A -> B+ (i.e. ``\verb+A+ implies \verb+B+'') in one text file, and you separately
prove \verb+B -> C+ in another. You might think
you have proven \verb+A -> C+.  However, when you merge the
two files together to get a proof of \verb+A->C+, the verification program might
come up with a \verb+Universe Inconsistency+ error
if the constraints on the universe
indices which are generated by the proofs of \verb+A->B+ and 
\verb+B->C+, while separately satisfiable, are not satisfiable by any common assignment. 

This phenomenon can be seen as a ``disturbance due
to measurement''. If we think of \verb+A+, \verb+B+ and \verb+C+ as
measuring posts along a lab bench (such as polarized glass panes,
say \verb+A+ is horizontal, \verb+B+ is at $45$ degrees and \verb+C+ is
vertical)
then the result of starting with a photon polarized as \verb+A+,
doing the proof of \verb+A-> B+
and then ``measuring'' (i.e. stopping the text file and
saying ``wow, we just proved \verb+B+''), then inputting the
photon polarized as \verb+B+ into the next proof 
of \verb+B -> C+, we get a photon polarized as \verb+C+ in the output;
whereas if we do the whole proof starting from \verb+A+ and
trying to get out \verb+C+, no proof-photon comes through! 

In \cite{stm} \verb}qua.v}, a simple example of this phenomenon 
is given using Russell's paradox, assuming the excluded middle and using \verb}eqT_rect} 
via Maggesi's proof
of \verb}JMeq_eq} \cite{Maggesi} \cite{Cuihtlauac}. As pointed out by H. Boom
\cite{boom} it is possible (although more complicated) to do the same thing without 
axioms other than the basic \verb+Coq+ system, using the Burali-Forti paradox
\cite{Girard} \cite{BarrasCoquandWerner}.

It is important to note, then, that as long as the
typical ambiguity mechanism is in place, it is not
allowable to prove a theorem in one file and 
import it as an axiom in another file.  
It is actually {\em required} that the proofs of
component pieces of the argument be imported by cut-and-paste (2.4).
Luckily the system has a \verb+Require+ command that does
this automatically without having  to copy the files.

Questions about ``universe polymorphism'' lead to some of the
most difficult notational problems. The following
discussion comes from work of Marco 
Maggesi, who is also working on solutions to the
problem possibly by implementing the thesis of Judica\"{e}l Courant \cite{Courant}.  
As was  pointed out in \cite{Feferman}, the basic example is the distinction between
small categories and big categories, in other words  assigning universe levels in the definition
of category.  Ideally, polymorphism should
be extendable to definitions such as this. In the
absence of such a mechanism, we seem to have to rewrite
\verb}Category_i} once for each level \verb+i+. This can lead to a nightmare
of notation, for example do we need to distinguish between
a \verb+comp_sm+ for composition in a small category, and
\verb+comp_bg+ for a big category? If
this seems annoying but doable, notice that we then
need \verb}Functor_ij} for functors from
a \verb}Category_i} to a \verb}Category_j}; and 
composing them, and so forth!   

One of my main motivations for going back to a purely set-theoretic
approach (cf 10.3, 10.4 below) is to try to avoid the above problems as much
as possible. However, 
it seems that there will always be a threshold between
the objects we are willing to manipulate and the bigger ones
we would prefer not to manipulate, and looking at categories of the
smaller objects inevitably bumps up against the bigger ones.
Until someone comes up with an adequate general
solution (perhaps the ideas in
\cite{Feferman} are relevant)  we are left with little choice but to live with
the situation.

\section{Exotic versus regular math}

The system described above applies to a wide variety
of situations, problems and logics. 

\numbering {\sc Intuitionism}: One of the main utilizations of this
type of system is to verify intuitionistic logic. In fact
it is difficult to imagine taking a complicated argument
and trying to insure, as a human reader, that the proof
{\em didn't} use the excluded middle in some hidden
place in the middle. In this sense, computer verification
seems crucially important and allows a much more serious 
approach to the field.  

\numbering {\sc Impredicative sort}: in previous versions of \verb+Coq+ the sort 
\verb+Set+ was declared as {\em impredicative},\footnote{In v8.0 the impredicativity
is turned off and can apparently be turned back on by an option.}
in other words if \verb+X:Set+ then the product type
\begin{verbatim}
Set -> X := (Y:Set)X
\end{verbatim} 
is classified as having type \verb+Set+.
Nonetheless
elements \verb+X:Set+ behave somewhat like sets: strong elimination is allowed, so that 
constructors in inductive sets are distinct. 
The intuition behind this situation is that the constructive
functions on a constructive set themselves form a constructive
set too. Thus these considerations are closely related to
the constructive or algorithmic aspects of mathematics. 
Hugo Herbelin
points out that this can lead to some unexpected and
fascinating things, such as a proof which extracts to 
an algorithm where nobody can understand how it
works \cite{Herbelin}. 

\numbering {\bf Caution:} the existence of an impredicative but
strongly eliminatory sort is somewhat ``orthogonal'' to 
standard
mathematics, in that it contradicts most 
choice axioms you can imagine, even the much weaker 
dependent
choice axioms \cite{Geuvers} \cite{ChicliPottierSimpson}.
Perhaps not all, though: David Nowak in \cite{Nowak} 
proposed doing the axiom of choice for the \verb}Ensembles} library,
which means looking at a function associating to any predicate 
\verb}P:X->Prop} another predicate \verb}(choice P):X->Prop}
with the property that if \verb}P} is nonempty then 
\verb}(choice P)} is a singleton. One conjectures that this
version is consistent even with an impredicative \verb}Set}.
Nowak's version is difficult to use in practice, and leads to a 
situation where the only objects which are manipulated are predicates.
This nullifies much of the advantage of using a type-theory framework.
So, for implementing classical mathematics 
we prefer working without the impredicativity. 

\numbering {\sc Logic in a topos}: recall from Moerdijk-MacLane
\cite{MoerdijkMacLane}
that one way of providing a model for certain intuitionistic
logical systems is by looking at logic in a topos other
than the topos of sets.  One can say this differently
(as they do in the book) by saying that any logical 
reasoning which is expressed solely in a limited intuitionist
formulation, applies internally to logic within any topos.
One can formulate as a conjecture, that this principle extends
to the full intuitionist logical system as implemented in \verb}Coq}:

\noindent
{\bf Conjecture}\, {\em 
The $\lambda - \Pi$ calculus
with universes and inductive objects as described above and
implemented in \verb}Coq}, applies to logic within a 
Grothendieck topos.}

It would be a generalization of the principles set
out in \cite{MoerdijkMacLane}. 
Some part of it may already be known, see Streicher \cite{Streicher}
as well as the references therein (and indeed the formulation of the question
appears to have been known to workers in Synthetic Domain Theory for some time
\cite{FioreRosolini} \cite{ASimpson}). 
It is not immediate, for
example, to see how to implement the sorts
\verb}Type_i}.  A {\em caveat} might also be in order because
recent versions of \verb}Coq} include some
inductive principles such as \verb}False_rect} or
\verb}eqT_rect} which might not be realizable in the internal
logic of a topos.  So, the above statement is really 
just intended to point out that it is an interesting 
research problem to see how far we can extend the internal-logic
principle from \cite{MoerdijkMacLane}.

As an example of how this type of thing might actually
lead to useful results for standard mathematics, one
can use the full flexibility of the \verb}Coq} implementation
of type theory to define internal cohomology within
the logic (this is done for $H^1$ in \cite{stm} \verb}h1.v}).
If this could be applied to a Grothendieck topos, it would
show that cohomology in the topos was an internal logical 
construction. Internal algebraic geometry inside the 
topos of sheaves on the big fpqc site takes on a particularly
nice form very close to the original Italian geometry,
for example affine $n$-space is quite literally just the
cartesian power of the ground field $k^n$. 

\numbering {\sc The Hodge conjecture}:  The statement of the
Hodge conjecture implies that the problem of saying whether a given topological cycle is algebraic,
is decidable. Same for the question of whether a given cycle is Hodge (i.e. vanishing of the Hodge
integrals).  This might be considered as a logical analogue of 
Deligne's observation that the Hodge conjecture implies that
Hodge cycles are absolute.  The pertinance of calculatory methods
showed up recently in \cite{Shimada}. It could be interesting to investigate further, for example does Deligne's proof
of \cite{abs} imply the decidability for abelian varieties? And in exactly which level of intuitionist
logic would the decidability hold? What are the consequences of this decidability given already by
the Lefschetz $(1,1)$-theorem? What happens internally in a topos?
Could one obtain a counterexample to the Hodge conjecture by showing
that in some cases the question is undecideable?

\begin{center}
--------------------
\end{center}

\numbering So there are undoubtedly lots of fun things to look at in these exotic reaches, but
they don't really have anything to do with our original goal of implementing,
as quickly as possible, a large amount of standard mathematics.
More generally, my experience is that
any ``brilliant gadget'' seems predestined
to be a bad idea. Boring as it may seem, we have
just to plod along in a mundane but systematic way. 

\section{Typed versus untyped style}

One of the central questions in the theory has been
whether to impose a strict type system, or whether to relax the
typing constraints. 
There are some very interesting threads of discussion about the issue
in the mailing lists, particularly \cite{qed} following the postings
\cite{Lamport} and \cite{Holmes}.  It is interesting to note that in
the first thread of discussion, Lamport (forwarded by Rudnicki \cite{Lamport}) offered
the viewpoint that the untyped style was better than the typed one.
He was promptly refuted by almost all the subsequent messages. A year
later when Holmes \cite{Holmes} asked what options people generally 
used to answer the question, almost everybody replied that the untyped style
was obviously the way to go. Unfortunately, this seems to have marked
the end of interesting discussion on the QED mailing list \cite{qed}. At 
about the same time the Coq system (where
typed language was systematic) started picking up speed: the announcement of 
Coq Version 5.10 is among the last group of messages in Volume 3 of \cite{qed}. 

The published version of \cite{Lamport} is \cite{LamportPaulson}, issued from
a combination of Lamport's original note and a rebuttal by L. Paulson. This discussion
and the discussion in J. Harrison \cite{Harrison} are excellent references
for the multiple facets of the question which we don't attempt further to reproduce
here. 

The distinction between
the typed and the untyped style is generally seen as related to the axiom of choice. In the absence of the
axiom of choice (in particular, in any sort of constructive environment) it
is much more natural to keep a lot of type information, looking at the values of functions only  
when they are well defined and make sense.
On the other hand, the axiom of choice provides a nice way to implement 
aspects characteristic of an untyped environment such as function evaluation which
is defined everywhere \cite{LamportPaulson}.  

Some authors make a distinction between choice and dependent choice, the latter stating
only that we can choose a function whose values lie in a family of single-element types.
By \cite{Geuvers} \cite{ChicliPottierSimpson} this weaker form poses the same problems
vis-a-vis a constructive axiom such as impredicativity of \verb}Set}. Furthermore
it generates a lot more proof obligations since one has to prove unicity whenever it is used.
For this reason we don't pay too much attention to the distinction, and 
look directly at the full axiom, even a strong Hilbertian form involving
an $\epsilon$-function.

\section{Options for doing set theory}

Sets are the most basic building blocks entering into the
definition of various mathematical objects. On the
other hand, as we have seen above, natural considerations
about general mathematical theories lead to type theory
as the framework for the system as a whole. The most critical
decision about document style is therefore the way the
notion of a set is encoded in the ambient type theory. 

We will review some of the most prominent options. These follow the contours
of the ``typed versus untyped'' discussion. 

\numbering 
{\bf Interpret sets as types}.  Benjamin Werner's paper
providing the justification for the $\lambda-\Pi$ calculus
with inductive definitions \cite{Werner} rests
on a {\em semantics}, i.e. assignation of meaning to the
formal symbols of the theory, where \verb}X:Type_i}
corresponds to a set in the $i$th Grothendieck universe
$U_i$, and functions such as constructed by the $\lambda$ operation
\verb}fun x:A =>( ... x ... )} are represented by sets which
are the graphs of the corresponding set-theoretic function.
We can use this semantics as our method for interpreting
sets, namely when we want to talk about a set we just take
any \verb}X:Type}.  In this option we make full use of the notion 
of dependent types, see C. McBride \cite{AltenkirchMcBride}. The following example shows how we would
define the type of groups. It uses the \verb}Record} construction
which is a special inductive definition corresponding to an
explicit list (of objects possibly depending on the previous
ones). 

{\scriptsize
\begin{verbatim}
Record Group : Type_(i+1) := {
elt:Type_i; 
op : elt -> elt -> elt; 
id : elt;
inv : elt -> elt;
assoc : forall x y z:elt,(op (op x y) z) = (op x (op y z));
left_id : forall x:elt, (op id x) = x;
right_id : forall x:elt, (op x id) = x;
left_inv : forall x:elt, (op (inv x) x) = id;
right_inv : forall x:elt, (op x (inv x)) = id
}.
\end{verbatim}
}

With this definition, if \verb}G:Group} then 
\verb}(elt G):Type_i} is the ``set'' of elements of \verb}G}.

The advantage of this viewpoint is that it is generally quite easy
to understand what is going on. You the uninitiated reader 
didn't have any trouble reading
the previous definition, and could by cut and paste 
give the analogous definition of, say, \verb}Ring}, right?

One of the drawbacks with this point of view is that if 
\verb}H:Group} too, and if for some reason we know
\verb}G=H}, we have problems interpreting elements of
\verb+G+ as being the same as elements of \verb+H+.
For example it is not legal to write 
{\scriptsize
\begin{verbatim}
forall x y z:(elt G), (op H (op H x y) z) = (op G x (op G y z))
\end{verbatim}
}
\noindent
(if rare at the outset, this sort of thing tends to come up in the middle
of proofs). 
The problem is alleviated by the recent inclusion in \verb}Coq}
of the ``transport'' function \verb}eqT_rect}, the utilisation of which is best 
accompanied by a systematic exploitation
of the \verb}JMeq} equality proposed by McBride \cite{AltenkirchMcBride}. This doesn't make
the problem go away entirely, though. 

This point of view adheres pretty closely to the ``typed'' side
of the discussion in the previous section. 
An appropriate choice axiom could allow us to extend function definitions to 
larger types, whenever a default is available. Unfortunately this implies a
distinction between empty and nonempty types. Some have made more exotic proposals
to include a default element $\perp$ in every type (this is based on work of
Scott \cite{Scott}, and was 
mentionned in the replies to \cite{Holmes}). I  haven't looked closely at it,
mostly for fear of a large number of proof obligations of the form $x\neq \perp$,
and because it obviates the advantage that we can easily understand what is going on. 
Also because the guy in the office next door talks about it so much. 

To close this subsection we also mention another more subtle drawback, which might
really be the main problem. There is a very slight shift in universe indexing, compared
with what we are used to in set theory. In the above \verb}Record}
if we want the underlying set \verb}(elt g)} to be
in the universe \verb}Type_i}, the type \verb}Group} is forced to be in the universe
\verb}Type_(i+1)}.  In general, a set at a given universe level will also have elements
at that level. Thus, {\em a priori}, the elements \verb}g:Group} are to be considered 
(even though this statement doesn't have any precise technical meaning) as being
in the $(i+1)$-st level. 
This contrasts with set theory where, if \verb}g} is a group whose underlying set
\verb}(elt g)} is in the $i$-th universe, then \verb}g} is represented by an
ordered uple which is also in the $i$-th universe.

\numbering 
{\bf Setoids}.  From a constructive point of view, the
notion of a quotient can be problematic because an
equivalence relation need not be decidable, so when
defining an object as a carrier modulo a relation,
the carrier might be constructive or computational but
the quotient less so. Another problem is that the classical
construction of a quotient as a set of equivalence classes
can pose problems in certain intuitionist frameworks 
(imagine doing that internally in a topos, for example).

 For these (or perhaps other) reasons a much-favored implementation of 
sets is
as {\em setoids} (cf the \verb}Setoid} libraries in 
\verb}Coq}). A setoid is a carrier type together with
an equivalence relation. We can think of them as quotient sets
without adding additional axioms. 
One major drawback is that there is a lot of extra
stuff to prove, basically that all operations defined on
the level of the carrier type are compatible with the
equivalence relation.  

For standard mathematics, we don't
really need the setoid viewpoint since we are willing to 
add in
all of the axiomatic baggage which makes the usual
theory of quotients go through.  
What is useful is to know that such structures exist, are natural, and
have been studied, because they occasionally show up uninvited. 

\numbering 
{\bf Introduce a type \verb}E} of sets}. 
In this viewpoint, the type-level is considered
as the meta-language of mathematical expressions, and
we introduce maybe some parameters
\begin{verbatim}
Parameter E : Type.
Parameter inc : E -> E -> Prop.
\end{verbatim}
Here \verb}E} is for ``Ensemble'' and \verb}inc} is
the elementhood relation. 
One can then go on to suppose
all of the usual axioms in set theory, and proceed from there.
For example the definition of the
subset relation becomes
\begin{verbatim}
Definition sub := [a,b:E](c:E)(inc c a) -> (inc c b).
\end{verbatim}

A development of set theory using this type of
viewpoint was done by Guillaume Alexandre in his Coq
contribution ``An axiomatisation of intuitionistic Zermelo-Fraenkel 
set theory'' (1995) \cite{Alexandre}.
This version is also quite analogous to the foundations of the MIZAR system \cite{Mizar},
and has been developped rather deeply in the ``ZF'' logic in the Isabelle
system \cite{Isabelle}. 

One can even {\em create} such an \verb}E} using an inductive definition
and a quotient, see Aczel \cite{Aczel} and Werner
\cite{Werner}:
\begin{verbatim}
Inductive Ens : Type := sup : (A:Type)(A->Ens)->Ens.
\end{verbatim}
This is a carrier type, bigger than what we actually think of as $Ens$ because
an equivalence relation of extensional equality is defined (thus, \verb}E} should actually be thought of as
a setoid and to get a type we would have to take the quotient).
Once we have the 
basic axioms, which in these constructions may
require supposing their counterpart axioms on the type-theory
side \cite{Werner},
the further development of the theory doesn't require 
knowing about the construction, so in this sense it is
conceptually easier just to take \verb}E} as a parameter.

In theories of this option, sets are 
terms of type \verb}E}, but not themselves types, and elementhood is just a propositional relation
on \verb}E}.  This is contrasted with the
type-theoretical interpretation (10.1) where the
relation \verb}x:X} is meta to the theory, not
actually accessible from within, and is decided by the computer. 

This would seem at first glance
to be much less powerful than version (10.1) where 
sets are types.
In a certain sense that's true, but in some other ways
it is actually more powerful. The fact that the elementhood
relation is a propositional function in the theory
allows us to manipulate it mathematically.
The problem referred to above about elements
of sets which are equal, goes away: if \verb}a=b} and
if \verb}x:E} then the propositions \verb}(inc x a)}
and \verb}(inc x b)} are equal since they are functions
of \verb+a=b+. Inside proofs, the fact that a certain
term has a certain type is replaced by the 
statement that a certain set is an element of another set;
in the first case the knowledge is required at the start
because otherwise the statement
will be rejected by the type-checker, but in the second case the proof can
be deferred into a separate branch of the proof tree.

If we further add an appropriately
strong version of the axiom of choice (basically 
Hilbert's $\epsilon$)
then we get an even better situation \cite{Lamport} \cite{LamportPaulson}
where function
evaluation can be defined as the operation which to any \verb}f,x:E}
assigns a choice of \verb}y:E} such that the pair
\verb}(pair x y)} is contained in \verb}f}, if such
a \verb}y} exists (or a default value if no such exists). 
In fact, if we write down all of the axioms in this
general spirit of everywhere-definedness, then  
everything that follows will be everywhere-defined too,
and we are no longer bothered with objects depending on proofs of propositions.

The main disadvantage of the purely set-theoretical approach is that we have
lost any possibility for the computer directly to compute 
on sets, because computation requires 
specific knowledge about the inductive structure of
the type containing a given term and our \verb}E}
no longer has any inductive structure. (This might
well be alleviated using the direct constructions of
\cite{Aczel} and \cite{Werner}, though, a thought-experiment
that starts off the next option.)

\numbering {\bf Combine 10.1 and 10.3}. 
In the Aczel-Werner type constructions 
\cite{Werner} one can construct a set i.e. term
\verb}x:E} whenever one has a type \verb}A:Type} together
with a realization function for its elements
\verb}A -> E}.  Imagine that we want to realize as a set
an inductive type say 
\begin{verbatim}
Inductive nat : Type :=  O : nat | S : nat -> nat.
\end{verbatim}
To do this we would need to construct 
an injection \verb}nat -> E}.
In this particular case it is natural to do so because
the construction 
$$
\emptyset , \;\; \{ \emptyset \} , \;\; 
\{ \emptyset , \{ \emptyset \} \} , \ldots 
$$
played an important historical role in the beginning.
However, if we try to do the same for, say, a free monoid,
few people ever really worry about the details of how to do
it. This requirement would therefore seem to qualify
as something that doesn't have any real mathematical
significance: we only care about the universal properties
of an inductive type, not about how it is actually realized.

In Werner's paper \cite{Werner}, he describes briefly
how to realize inductive
types as sets, and consequently how to realize every type as
a set. So, let's just abstract out this process
and suppose, with a parameter, that we know how to
do it but we don't care what it is. We get a situation
where, combining (10.1) and (10.3), sets are types,
but also the elements of sets are realized as sets
themselves---which after all is the bedrock assumption of classical
set theory. 

This is the viewpoint which is developped in  \cite{stm}.  We start by putting
\begin{verbatim}
Definition E:= Type.
Parameter R : (x:E)x->E.
Definition inc := [a,x:E](exists [b:x](R b) = a).
\end{verbatim}
Thus, we identify the notion of set with the notion of type. 
The parameter
\verb}R} is the realization function which says that the elements of
a set/type are themselves sets. The elementhood relation \verb}inc} 
is defined using \verb}R} by
saying that $a\in x$ iff there exists a term $b$ of type $x$ whose
realization is $a$. 

We assume a certain number of axioms such as choice
(in the form of a Hilber $\epsilon$-function) and
replacement.  We can then go on to develop set theory for 
\verb}E} just as would be done in option (10.3) above.
On the other hand, we can calculate with inductive types
and other type-theory features, useful among other things 
to establish a notational system.

\section{Inheritance}

Consider the notion of Lie group. This fits into at least
4 distinct (plus one overlapping) inheritance threads as shown in Figure 1.
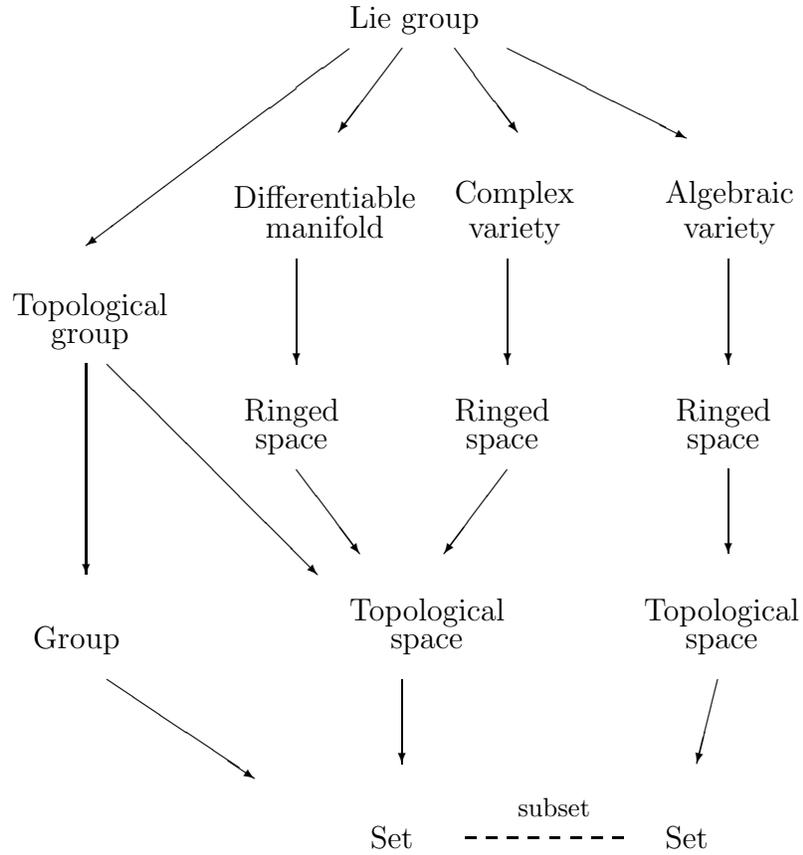
\begin{figure}
\setlength{\unitlength}{1.4cm}
\begin{picture}(6,9)(.7,-.3)
\put(3,8){\shortstack{Lie group}}
\put(3,7.8){\vector(-4,-3){2.5}}
\put(3.5,7.8){\vector(-3,-4){.6}}
\put(4,7.8){\vector(3,-4){.6}}
\put(4.5,7.8){\vector(2,-1){1.7}}
\put(1.9,6){\shortstack{Differentiable\\ manifold}}
\put(4,6){\shortstack{Complex\\ variety}}
\put(6,6){\shortstack{Algebraic\\ variety}}
\put(-.2,5){\shortstack{Topological\\ group}}
\put(2.5,5.8){\vector(0,-1){1}}
\put(4.5,5.8){\vector(0,-1){1}}
\put(6.6,5.8){\vector(0,-1){1}}
\put(.5,4.8){\vector(0,-1){2}}
\put(.7,4.8){\vector(1,-1){2}}
\put(2,4){\shortstack{Ringed\\ space}}
\put(4,4){\shortstack{Ringed\\ space}}
\put(6.1,4){\shortstack{Ringed\\ space}}
\put(2.5,3.8){\vector(3,-4){.6}}
\put(4.5,3.8){\vector(-3,-4){.6}}
\put(6.6,3.8){\vector(0,-1){.8}}
\put(0,2.1){\shortstack{Group}}
\put(3,2.1){\shortstack{Topological\\ space}}
\put(5.8,2.1){\shortstack{Topological\\ space}}
\put(.7,1.8){\vector(3,-2){1.4}}
\put(6.5,1.8){\vector(-1,-4){.2}}
\put(3.5,1.8){\vector(0,-1){.8}}
\put(3.2,.2){\shortstack{Set}}
\put(6,.2){\shortstack{Set}}
\multiput(4.1,.3)(.2,0){8}{\line(1,0){.1}}
\put(4.6,.5){\shortstack{{\footnotesize subset}}}
\end{picture}
\caption{Inheritance for the notion of Lie group}
\end{figure}
The three rightmost threads are formally similar,
however the ringed spaces in question are different,
and the topological spaces for the two middle
threads are the same but different from the 
(Zariski) topological space in the ``algebraic variety''
thread. One could identify other inheritance strings,
and note that there are intricate relationships
between them, for example the structure
sheaf of the ringed space underlying
the complex variety is a subsheaf of the 
structure sheaf of the ringed space underlying the
differentiable manifold, and the set underlying the usual topological 
space in the first threads
is the subset of closed points in the set underlying the Zariski
topological space. 

If $G$ is
a Lie group, informally we tend to use the letter ``G'' for
any and all of the above things and more.

On the other hand, the 
details of these inheritance relations would
have precisely to be spelled out in
any computer verification system.

Letting the same
letter ``G'' mean any of an amorphous mass of notions,
represents a huge collection of notational shortcuts.
While clearly impossible to keep all of them,  
it should be possible
to keep some. This could lead to some difficult
choices: how to decide which parts to keep and which parts
to relegate to a more complicated notation. And in any case,
how do we minimize the amount of additional notation
necessary to make these shortcuts precise?  

The example of a Lie group is deliberately complex in
order to illustrate the problem. However, it comes up
to varying degrees in pretty much all of modern 
mathematical practice.

\section{Perspectives}

We now outline a sketch of what some part of the further development of machine-verified
mathematics might look like.  This could start with any of a number of 
basic supports such as \cite{Isabelle}, \cite{Mizar}, \cite{Coq}, \cite{Metamath}, \ldots , some of which
are already quite far along.

Attached to the original preprint version of the present paper was the start of a development along the lines
of option (10.4) above. The idea of notative sets as being functions on strings, first brought up in 
\cite{LamportPaulson}, was implemented. We get a primitive kind of inheritance, in which different structures corresponding
to different strings can coexist in the same object, and the functions extracting these operations from
a given object are uniform throughout. Thus for example the theory of groups and the theory
of topological spaces would both apply immediately to topological groups. 
The explanatory section and source code have now been spun off as a separate
preprint \cite{stm}. 

The proof development accompanying \cite{stm} now treats the notions of ordinal and cardinal. The next part 
which needs to be addressed would contain the
notion of finiteness, various facts about finite sets (notably that any automorphism of a finite set
decomposes as a product of transpositions); also various notions of collection or list of objects; etc. 
The theory of transpositions should help for defining multi-term summation in an abelian monoid.

Most objects are constructed with a single common variable, denoted \verb}(U a)}, thought of as the 
``underlying set'' of \verb}a}.  Morphisms between such objects can be considered as some of the maps between underlying sets.
This makes it possible to develop a  certain amount of categorical machinery
at once for all these types of objects.

With this in mind, we can define monoids, groups, rings, algebras, modules; categories, eventually with 
Grothendieck topologies; presheaves (and sheaves if the category has a g.t.); ordered sets; 
topological spaces also eventually
with other structures so that one could define topological groups and so on.  

The counting system has to be extended to the numbering system.  We can import  much of the beautiful
development of Peano arithmetic already in the \verb}Coq} system, as well as the \verb}ZArith} library containing a
logarithmically fast implementation of integers.  Then we need to define the rationals as a localization of
the integers (and for this we might as well call upon the general notion of localization of a ring).  Then we
can define the reals as the set of Dedekind cuts of the rationals.  It is not clear whether it is better to
import the real number results from the existing \verb}Coq} library, or to redo a substantial portion, since our approach
using things like the axiom of choice, and an untyped point of view, might give substantial simplification here
(as opposed to the case of elementary arithmetic where we wouldn't expect any such benefit).  

Once we have the reals (and therefore the complex numbers too) we need to establish the basic facts about the
topology on these spaces, and along the way note that this opens up the definition of metric space, so we can
do the basic facts about metric spaces.  This could eventually lead to stuff for analysis such as 
Hilbert or Banach spaces, and all of the basic theorems of analysis.  This is not necessarily urgent although
it would be required for the implementation of Hodge theory. 

A parallel project is to develop abstract commutative and linear algebra as well as group theory and field theory.
For large parts, these theories don't need the real numbers.  As pointed out above, the notion of localization
would be convenient to have for the construction of the rational numbers.  In general, it would be good to start
by developing the basic commutative algebra such as is assumed in the beginning of Hartshorne.  Once we have
a reasonable control over the basic notions of ring, algebra, module etc. then much of commutative algebra can
be done without any new notational worries.  

Linear algebra could lead to a rather extensive development of the theory of matrices, matrix groups, and
group representations. 

Again in parallel, we need the notions of presheaf and sheaf over a topological space.  Here it seems
economical to take the abstract point of view and develop directly the notion of sheaf on a Grothendieck site.
Along the way the notion of presheaf on a category yields the notion of simplicial set which forms the basis
for homotopy theory. The notion which is directly useful for starting out is that of sheaf on a topological space,
so we need a construction going between a space and its corresponding Grothendieck site of open sets.  
It is  unclear whether we would prefer using the full abstract notion of point for a site (or maybe a topos?---but
this may lead to universe problems)
or whether it is better just to look at the germ for a presheaf on a topological space.  In any case, here some
{\em ad hoc} argumentation is necessary, because we would like the germ of a presheaf to be endowed with the
same type of structure as the presheaf; the variety of different types of structures we want to preserve means that it is probably
more difficult to come up with a single general theory than just to do each one separately.  

Up to here, most of what is listed above has already been treated in one way or another, in at least
some user-contributions for at least some of the systems listed in the bibliography. 
Left open is the problem of coalescing everything together. 

Once we have these various theories, we can put them together to obtain the notion of locally ringed space modelled
on a certain type of model space. Depending on which model spaces we look at, we can obtain the notions of
differentiable manifold, complex analytic manifold, or algebraic variety (scheme). 
Laurent Chicli has done a
formalization up to the notion of (affine) scheme in his thesis \cite{Chicli}, using setoids (option (10.2) above).

At this point, perhaps one of
the first things to do would be the comparison constructions going from
schemes (over ${\bf C}$) to complex manifolds to differentiable manifolds.  We should then be able to construct 
projective space, projective subvarieties, various differentiable manifolds, and so on.  
Another thing to do will be 
exterior algebra and differential calculus
over a manifold.

It would be interesting to try things out by applying this to the case of Lie groups, and seeing how far
into the theory one can go. 

Integration is probably the next place where we will run into a major
notational problem.  This is because the verb ``to integrate'' spans a highly diverse and 
notationally wide-ranging collection of ideas: there is simple Riemann (or Lebesgue) integration on the
real line (already implemented in \cite{ftc}), 
but also more general integration of a measure, or integration of a differential form over various
things such as a manifold, or a smooth chain; or finally integration of currents and the like. These 
different notions have a very high degree of overlap, however they are not quite the same, and
they require differing notations and even differing levels of notational complexity.  

With the theory of integration would come the possibility of getting at the fundamental properties of analysis
on manifolds. At this point we will need the theory of Hilbert and Banach spaces (already implemented in \cite{Metamath}),
and we should be able to
do the basic theorems such as Stokes' theorem, the Poincar\'e lemma and de Rham's theorem, Hartogs' theorem
for complex manifolds, the mean value theorem for harmonic functions, which gives
the maximum principle and unicity of analytic continuation,
the representation (and regularity) of cohomology by harmonic forms, 
leading up to the Hodge decomposition
for the cohomology of a complex K\"{a}hler manifold (to take an example from my own research interests). 
At this stage we will also want to have the basic notions of homotopy theory available.

Once all of the above is well in hand, we would be able to state and prove a first theorem which is a
really nontrivial application of this wide swath of mathematics: the fact that the first (or any odd)
Betti numbers of complex algebraic varieties are even.  

We should then be able to start looking at more advanced things going
toward current research problems. 
One area where we can expect important gains is homological
or homotopical algebra, for example J. Lipman points out that it would be good to 
formalize the computations in duality theory such as \cite{ATJLLipman}.  
More geometrical directions could include for example the construction of moduli spaces of vector bundles,
connexions and so forth. Or the investigation of existence problems for connexions over higher
dimensional varieties which could at first be viewed just as matrix problems in linear algebra.
 
Some day in the not-too-distant future we will have a computer-verified proof that the curvature of the
classifying space for variations of Hodge structure in the horizontal directions is negative, 
and all that it entails.

 \end{document}